\documentclass[letterpaper,12pt]{amsart}
\usepackage[notref,notcite]{}
\advance\oddsidemargin by-0.5in
\advance\evensidemargin by-0.5in
\advance\textwidth by 1in

\newcommand{\gateaux}{G\^ateaux}
\newcommand{\frechet}{Fr\'echet}

\newcommand{\R}{\mathbb{R}}
\newcommand{\real}{\mathbb{R}}

\newcommand{\e}{\varepsilon}

\newcommand{\ball}[2]{B_{#2}(#1)}
\newcommand{\ballcl}[2]{\overline{B}_{#2}(#1)}
\newcommand{\q}{q}
\newcommand{\T}{T}

\newcommand{\setseq}{\mathfrak S}

\renewcommand{\phi}{\varphi}
\newtheorem{theorem}[subsection]{Theorem}

\newtheorem{lemma}[subsection]{Lemma}
\newtheorem{definition}[subsection]{Definition}

\numberwithin{equation}{section}

\newcommand{\bib}{
}

\begin{document}
\title[Compact Universal Differentiability Set]{A Compact Universal Differentiability Set with Hausdorff Dimension One}
\author{Michael Dor\'e}
\address{Mathematics Institute\\ 
University of Warwick\\
Coventry CV4 7AL\\
United Kingdom}
\email{Michael.J.Dore@gmail.com}
\author{Olga Maleva}
\address{School of Mathematics\\
University of Birmingham\\
Edgbaston\\
Birmingham B15 2TT\\
United Kingdom}
\email{O.Maleva@bham.ac.uk}
\begin{abstract}
We give a short proof that any non-zero Euclidean space has a compact subset of Hausdorff dimension one that contains a differentiability point of every real-valued Lipschitz function defined on the space.
\end{abstract}
\advance\baselineskip by3pt

\maketitle
\section{Introduction}
\subsection{}
It is well known that if $f:\R\to\R$ is a Lipschitz
function, then $f$ is differentiable almost everywhere.
Zahorski, \cite{Z}, gives a full characterisation of the possible
sets of points of non-differentiability of a real-valued Lipschitz function defined on $\R$. In particular, it
follows that for any Lebesgue null set
$E\subseteq \R$  there exists a Lipschitz function
$f:\R\to\R$ that is non-differentiable
at every point of $E$.

Turning to higher dimensions, we may still conclude that real-valued functions defined on a finite dimensional Euclidean space are differentiable almost everywhere; this is Rademacher's theorem. However the converse implication no longer holds; in any Euclidean space of dimension at least two, there are sets of measure zero on which every real-valued Lipschitz function, defined on the space, is somewhere differentiable. Examples of sets satisfying this property --- which we hitherto refer to as
the universal differentiability property
--- were constructed by Preiss in
\cite{P} and the authors of the present paper in \cite{DM}.

It is proved in \cite{P} that if 
$E\subseteq\real^d$ is any $G_\delta$ set, i.e. an intersection of a 
countable family of open sets, such that $E$ contains every line segment 
that passes through any two points of some dense subset $S \subseteq \R^d$,
then $E$ has the universal differentiability property. One can check that 
not only can such a set $E$ be taken to be null in $\R^d$ for $d \geq 2$ 
but that such a set may also be chosen to have Hausdorff dimension one; 
see Lemma~\ref{hauslemma}. However it is also clear that the closure of 
such a set $E$ is the whole space $\R^d$.

In \cite{DM} we show, on the other hand, that it is possible to
find a compact and null subset $E \subseteq \R^d$ with the
universal property, for $d \geq 2$. An example for $d = 2$ is
given by a generalisation of the Menger-Sierpinski carpet;
see also \cite{IM}. More precisely we choose odd integers $N_i > 1$
such that $N_i\to\infty$ and $\sum N_i^{-2}=\infty$. At the first step, 
we divide the unit square $[0,1]^2$
into $N_1^2$
equal squares and remove the central square in the division. Then on each
subsequent step we divide the remaining squares into $N_i^2$ equal squares and remove the central square. The example is given by the set that remains. For $d > 2$ we take the Cartesian product of this two-dimensional set with $\R^{d-2}$. See \cite{DM} for more details.

The drawback of the example in \cite{DM} is, however, that its Hausdorff dimension is equal to $d$, the dimension of the underlying space.

In the present paper
we construct a compact
subset of $\R^d$ with the universal differentiability
property such that its Hausdorff dimension is equal to one,
thus making the set small both in terms of its closure
and in terms of its Hausdorff dimension.

For Lipschitz mappings to spaces of dimension larger than
one there are fewer positive results.
For $d\ge 3$, it is proved in \cite{PH}
that there exists a Lebesgue null
set $E$ in $\R^d$
such that every Lipschitz mapping from $\R^d$ to $\R^{d-1}$ has a
point of
$\e$-differentiability in that set for all $\e > 0$; see the subsequent section for a definition.
In fact one may take $E$ to be the union of all ``rational
hyperplanes'' in $\R^d$, so that the Hausdorff dimension of $E$ is equal to $d-1$.
However $\e$-differentiability is weaker than
differentiability.

\subsection{}
Recall for a pair of real Banach spaces $X,Y$
a function $f:X\to Y$ is called Lipschitz if there is a constant
$L \geq 0$ such that
$$\|f(x')-f(x)\|_Y\le L\|x'-x\|_X$$
for any
$x,x'\in X$. The smallest such $L$ is denoted as
$\text{Lip}(f)$.

We say that the function $f:X\to Y$ has a directional derivative
at $x\in X$ in the direction $e\in X$ if the limit
\begin{equation}\label{eq.dirder}
\lim_{t\to0}\frac{f(x+te)-f(x)}{t}
\end{equation}
exists. We then denote the limit \eqref{eq.dirder}
as $f'(x,e)$. We say $f$ is \gateaux{} differentiable
at $x\in X$ if $f'(x,e)$ exists for every $e\in X$ and
$T(e):=f'(x,e)$ defines a bounded linear operator $T:X\to Y$.

If $f$ is \gateaux{} differentiable at $x$ and the limit
\begin{equation}\label{eq.fre}
\lim_{h\to0}\frac{f(x+h)-f(x)-T(h)}{\|h\|}
\end{equation}
is equal to $0$ --- or equivalently, the
convergence in \eqref{eq.dirder}  is uniform for $e$ in the unit sphere of $X$ --- then we say that $f$ is \frechet{} differentiable
at $x$ and denote the operator $T$ as $f'(x)$.

The condition \eqref{eq.fre} can be rewritten as follows.
We require that there exists a bounded linear operator $f'(x)\colon X\to Y$ such that
for any $\e>0$ there exists $\delta>0$ such that for any
$h\in X$ with $\|h\|<\delta$ we have
$$
\|f(x+h)-f(x)-f'(x)(h)\|\le\e\|h\|.
$$
If, on the other hand, we only know the existence of such an operator for some fixed
$\e$, we say that $f$ is $\e$-\frechet{} differentiable at $x$.

We refer the reader to \cite{JLPS,LP} where the notion
of $\e$-\frechet{} differentiability is studied in relation to Lipschitz mappings, with the
emphasis on the infinite dimensional case.
In general, \frechet{} differentiability is a strictly stronger property than \gateaux{}
differentiability. However the two notions coincide for Lipschitz functions defined
on a finite dimensional space; see \cite{BL}. Hence, in this case, we may simply refer to differentiability, without any ambiguity.

The Hausdorff dimension of a set $E$,
a subset of a metric space, is defined in the following
way. For $r,\delta>0$, let
$$
\mathcal  H_\delta^r(E)=\inf \left\{
\sum_{i=1}^\infty[\textrm{diam}(S_i)]^r \text{ over all }S_i \text{ with }
E\subseteq \bigcup_{i=1}^\infty S_i
\text{ and }
\text{diam}(S_i)\le\delta \right\},
$$
and
\begin{equation}
\label{deltaless}
\mathcal H^r(E)=\lim_{\delta\to0^+}\mathcal  H_\delta^r(E);
\end{equation}
as $\mathcal H_{\delta}^r(E)$ is a decreasing function of
$\delta > 0$, this limit exists in $[0,\infty]$.
The number
\begin{equation}
\label{defofH}
\text{dim}_\mathcal H(E)
:=\inf\{r>0 \text{ with } \mathcal H^r(E)=0\}
\end{equation}
is called the Hausdorff dimension of $E$.
It is easy to see that the Hausdorff dimension is
a monotone set function with respect to inclusion
and
if $f \colon X\to Y$ is a Lipschitz function
then the Hausdorff dimension of $f(E)$ does not exceed the
Hausdorff dimension of $E$, for every $E\subseteq X$.
See \cite{Mat} for a discussion of the properties
of Hausdorff dimension.

From this point, we shall work in $\R^d$ where $d \geq 1$. We fix some notation for the rest of the paper. We let $\|\cdot\|$ be the Euclidean distance on $\R^d$ and
denote an open ball centered at $a \in \R^d$ of
radius $r$ by $\ball{a}{r}$ and a closed ball by
$\ballcl{a}{r}$. Further, for any $S\subseteq\real^d$ and $r\ge0$ we
let 
\begin{align*}
\ball{S}{r}&:=
\{x\in\real^d\text{ such that }\inf_{y\in S}\|x-y\|<r\}=
\bigcup_{y\in S}\ball{y}{r},\\
\ballcl{S}{r}&:=\{x\in\real^d\text{ such that }
\inf_{y\in S}\|x-y\|\le r\}
\end{align*}
denote
the open and closed $r$-neighbourhoods of $S$ respectively.

We require the following simple observation.

\begin{lemma}
\label{hauslemma}
If $L \subseteq \R^d$ is a countable union of line
segments then there exists a $G_{\delta}$ set $O \subseteq \R^d$
with $L \subseteq O$ such that the Hausdorff dimension of $O$
is equal to one.
\end{lemma}

\begin{proof}
It clearly suffices to show that the Hausdorff dimension of $O$ can be taken to be less than or equal to one.

Note that
if $I$ is a line segment of length at most $1$ in
a Banach space $Y$, $r,\delta > 0$ and $k \geq 4/\delta$ is a
positive integer then
\begin{equation}
\label{labnick}
\mathcal H_\delta^r(\ball{I}{1/k}) \leq k\cdot\left(\frac{4}{k}\right)^r = 4^r \cdot k^{-(r-1)},
\end{equation}
as we may cover $\ball{I}{1/k}$ with $k$ open balls whose radii are equal to $2/k$, i.e. with diameters $4/k \leq \delta$.

Now let $L \subseteq \R^d$ be a countable union of line segments. One may write
$L = \bigcup_{m \geq 1} L_m$,
where each $L_m$ is a line segment of length at most $1$.
Let
$$
O_n = \bigcup_{m=1}^{\infty}\ball{L_m}{1/2^{m+n}}
\qquad
\text{and}
\qquad
O = \bigcap_{n=1}^{\infty}O_n.$$

Note that $O$ is a $G_{\delta}$ subset of $\R^d$, containing $L$. To verify that the Hausdorff dimension of $O$ is no greater than one, it suffices, by \eqref{deltaless} and \eqref{defofH}, to show that
$\mathcal H_\delta^r(O) = 0$
for every $\delta > 0$ and $r > 1$.
If $2^{n+1} \geq 4/\delta$, then
\begin{multline*}
\mathcal H_\delta^r(O)
\leq \mathcal H_\delta^r(O_n)
\leq \sum_{m=1}^{\infty}\mathcal H_\delta^r(\ball{L_m}{1/2^{m+n}})
 \leq \sum_{m=1}^{\infty}4^{r} \cdot 2^{-(m+n)(r-1)}\\
 = 4^r \cdot \frac{2^{-(n+1)(r-1)}}{1-2^{-(r-1)}},
\end{multline*}
using the countable subadditivity of $\mathcal H_\delta^r$,
\eqref{labnick} and $r > 1$.
Letting $n \to \infty$ we obtain
$\mathcal H_\delta^r(O) = 0$
as required.
\end{proof}

\subsection{}
We have already mentioned that by \cite{P}
any $G_{\delta}$ set $O$ that contains every line segment passing through two points of a dense subset $R \subseteq \R^d$ has the universal differentiability property and that if $R$ is countable, $O$ may be taken to have Hausdorff dimension one by Lemma~\ref{hauslemma}. Our strategy then is to construct a closed and bounded subset of such a set $O$ that still has the universal differentiability property; this will give our example of a compact universal differentiability set with Hausdorff dimension one.

The basic idea of the construction is as follows.
We write $O = \bigcap_{k \geq 1} O_k$ where $O_k$ are open subsets of $\R^d$ with $O_{k+1} \subseteq O_k$ for each $k \geq 1$. Then for each $k \geq 1$ we construct a family of closed subsets $M_k(\lambda)_{\lambda \in [0,1]}$ of $O_k$ with the property that $M_k(\lambda) \subseteq M_k(\lambda')$ for $\lambda \leq \lambda'$. Taking the intersection $T_{\lambda} := \bigcap_{k \geq 1} M_k(\lambda)$ we note that each $T_{\lambda}$ is a closed subset of $O$ and that $T_{\lambda} \subseteq T_{\lambda'}$ for $\lambda \leq \lambda'$. We then prove, using the details of the construction, that the family $(T_{\lambda})_{\lambda \in [0,1]}$ contains, in a certain sense, a large amount of line segments connecting two points in the dense set $R$.

Next, by quoting Theorem~\ref{th.incr}, which is Theorem~3.1 in
\cite{DM}, we show that given a Lipschitz function $f\colon
\R^d\to\R$ we can find a point $x \in T_{\lambda}$ for some
$\lambda < 1$
and a direction $e \in S^{d-1}$, the unit sphere of $\R^{d}$,
such that the directional derivative $f'(x,e)$
is almost locally maximal: if $\e > 0$ and
$x'\in\T_{\lambda'}$ is close to $x$, with $\lambda'$ sufficiently
close to $\lambda$, and $e' \in S^{d-1}$ is a direction such that
$(x',e')$ satisfies certain additional constraints, then
$f'(x',e') < f'(x,e)+\e$.

Finally we then prove $f$ is differentiable at $x$ with derivative
$$
f'(u) = f'(x,e)\langle u,e \rangle
$$
using Lemma~\ref{lemdiff}, which we quote from
\cite[Lemma~4.3]{DM}. This last step makes essential use of the fact that the family $(T_{\lambda})_{\lambda \in [0,1]}$ contains sufficiently many line segments.

We finish this section by noting that the Hausdorff dimension of one is optimal.
\begin{lemma}\label{lem.low}
If $d\ge1$  and
$E\subseteq \R^d$ is a universal differentiability set, then
the Hausdorff dimension of $E$ is at least one.
\end{lemma}

\begin{proof}
Assume $E$ has Hausdorff dimension strictly less than $1$.
Let $v$ be any unit vector in $\R^d$ and set
$$\phi(x)=\langle x,v\rangle.$$
Since $\phi \colon \R^d \to \R$ is Lipschitz, we conclude
$$\dim_{\mathcal H}(\phi(E))\leq\dim_{\mathcal H}(E)<1;$$
in particular  $\phi(E)\subseteq \R$ has Lebesgue measure
$0$. Hence there exists a Lipschitz function $g:\R\to\R$ that is
non-differentiable at every $x\in\phi(E)$. Then
$f := g \circ \phi$ defines a Lipschitz function from
$\real^d$ to $\R$
that is not differentiable at every $x\in E$, as the directional derivative $f'(x,v)$ does not exist for $x \in E$.
\end{proof}

\section{Construction}
\label{sec.main}

Let $R$ be a countable dense subset of $\ball{0}{1}$ and for each $\e > 0$ let $R(\e)$ be a finite subset of $R$ such that for every $x \in \ball{0}{1}$ there exists $r \in R(\e)$ with $\|r-x\| < \e$.

By Lemma~\ref{hauslemma} we may pick a $G_{\delta}$ set $O \subseteq \ball{0}{1}$ of Hausdorff dimension one such that $[r,s] \subseteq O$ for every $r,s \in R$. We write $O = \bigcap_{k=1}^{\infty}O_k$ where $O_k$ are open
subsets of $\real^d$ with $O_{k+1} \subseteq O_{k}$ for each
$k \geq 1$.
\begin{definition}
\label{defR}
For $k \geq 0$ we define compact sets $R_k \subseteq O$ and $w_k > 0$ as follows.
Let $w_0 > 0$ and $R_0$ be any compact subset of $R$; for example $w_0 = 1$ and $R_0 = \emptyset$.
For each $k \geq 1$ we let
\begin{itemize}
\item
$R_{k} = \bigcup \{[r,s] \text{ where }r,s \in R(w_{k-1}/k)\}
\cup R_{k-1}$,
\item
$w_{k} \in (0,w_{k-1}/2)$ be such that ${\ballcl{R_k}{w_{k}}} \subseteq O_{k}$.
\end{itemize}
\end{definition}
Since $[r,s] \subseteq O$ for every $r,s \in R(w_{k-1}/k) \subseteq R$ and $R(w_{k-1}/k)$ is finite, the set $R_k$ defined above is compact and is a subset of $O$. Then we may pick $w_k > 0$ as given above because $R_k \subseteq O \subseteq O_k$, $R_k$ is compact and $O_k$ is open.
Note that $(w_n)$ is a decreasing sequence that tends to zero.
\begin{definition}
\label{defM}
If $k \geq 1$ and $\lambda \in [0,1]$ we set
\begin{equation}
M_k(\lambda) = \bigcup_{k \leq n \leq (1+\lambda)k}
{\ballcl{R_n}{\lambda w_{n}}}.
\end{equation}
\end{definition}

We note that for each $k \geq 1$ and $\lambda \in [0,1]$, the set $M_k(\lambda)$ is a finite union of closed sets, so closed. As ${\ballcl{R_n}{\lambda w_{n}}} \subseteq O_n$ for any $n \geq 1$ we have, for $k \geq 1$,
\begin{equation}
\label{MinO}
M_k(\lambda) \subseteq \bigcup_{k \leq n \leq (1+\lambda)k}O_n = O_k
\end{equation}
as $O_{n+1} \subseteq O_{n}$ for all $n \geq 1$.

\begin{definition}
\label{defT}
Given $\lambda \in [0,1]$ we set
\begin{equation}
\label{defS}
\T_{\lambda} = \bigcap_{k=1}^{\infty} M_k(\lambda).
\end{equation}
\end{definition}

Note from \eqref{MinO} that $\T_{\lambda} \subseteq O_k$ for every $k \geq 1$ so that $\T_{\lambda} \subseteq O$ for every $\lambda\in[0,1]$,
and the set $O$ is bounded and has Hausdorff dimension one. 
 Further, as $\T_{\lambda}$ is an intersection of closed sets, it is closed. We conclude that for every $\lambda\in[0,1]$, the set $\T_\lambda$ is a compact subset of $\R^d$ of Hausdorff
dimension at most one.
Finally we note that if $\lambda_1 \leq \lambda_2$ then we have $\T_{\lambda_1} \subseteq \T_{\lambda_2}$
and if $\lambda=0$, Definition~\ref{defM} implies 
$M_k(0)=R_k$, so that by Definition~\ref{defT} we have 
$\T_0=\bigcap_{k=1}^\infty M_k(0)=R_1$.

\begin{lemma}
\label{shift}
Suppose that $0 \leq \lambda < \lambda+\psi \leq 1$ and
$x \in M_k(\lambda)$ where $k\geq 1$. If
$$
0 < \Delta \leq \psi w_n
$$
for all $n \leq (1+\lambda)k$ then we have
$$
\ball{x}{\Delta} \subseteq M_k(\lambda+\psi).
$$
\end{lemma}

\begin{proof}
Using Definition~\ref{defM} we may find $n$ such that
$k \leq n \leq (1+\lambda)k$ and
$x \in {\ballcl{R_n}{\lambda w_n}}$.
Noting that $\Delta \leq \psi w_n$ we have
$$
\ball{x}{\Delta} \subseteq
{\ballcl{R_n}{\lambda w_n+\Delta}}
 \subseteq
{\ballcl{R_n}{(\lambda+\psi) w_n}}
 \subseteq M_k(\lambda+\psi)
 $$
using Definition~\ref{defM} once more and $k \leq n \leq (1+\lambda+\psi)k$.
\end{proof}

\begin{lemma}
\label{crit}
Suppose that $0 \leq \lambda < \lambda+\psi \leq 1$. 
If $\eta \in (0,1)$, $\Delta > 0$ and $k \geq 1/(\psi \eta)$ satisfy
$\Delta > \psi w_n$ for some $n \leq (1+\lambda)k$ then there exists $\alpha \in (0,\eta \Delta)$ such that
$[r,s] \subseteq M_l(\lambda+\psi)$
for every $r,s \in R(\alpha)$ and $l \geq k$.
\end{lemma}

\begin{proof}
As the sequence $w_n$ is decreasing we may assume
$k\le n \le (1+\lambda)k$
so that using $k\geq1/(\psi \eta)$, we get
$$
\alpha := 
\frac{w_n}{n+1} < \frac{\Delta/\psi}{1/(\psi \eta)} = \eta \Delta.
$$
Let $l \geq k$ and choose any $r,s \in R(\alpha)$ so that $[r,s] \subseteq R_{n+1}$ by Definition~\ref{defR}. Note that as
$n\le(1+\lambda)k$ and $\psi k \geq 1/\eta \geq 1$,
$$
(1+\lambda+\psi)l
\geq (1+\lambda)k + \psi k
 \geq n+1.
$$
Thus we may pick $m \geq n+1$ with $l \leq m \leq (1+\lambda+\psi)l$. Then $[r,s] \subseteq R_{n+1} \subseteq R_m$. Hence $[r,s] \subseteq M_{l}(\lambda+\psi)$ by Definition~\ref{defM}.
\end{proof}

\begin{lemma}
\label{lem.main}
For each $\eta,\psi > 0$ there exists
$$
\Delta_0 = \Delta_0(\eta,\psi) > 0
$$
such that if $\Delta \in (0,\Delta_0)$,
$0 \leq \lambda < \lambda+\psi \leq 1$ and $x \in \T_{\lambda}$ there exists $\alpha \in (0,\eta \Delta)$ such that for every $r,s \in R(\alpha) \cap \ball{x}{\Delta}$ we have $[r,s] \subseteq \T_{\lambda+\psi}$.
\end{lemma}

\begin{proof}
Pick $\Delta_0 > 0$ with $\Delta_0 < \psi w_n$ for every
$n \leq 2/(\psi \eta)$. Now suppose that
$\Delta \in (0,\Delta_0)$. Pick a minimal $k \geq 1$ such that
$\Delta > \psi w_n$ for some $n \leq (1+\lambda)k$. Note that
$(1+\lambda)k > 2/(\psi \eta)$ so that $k > 1/(\psi \eta)$. Thus by Lemma~\ref{crit} we can find
$\alpha \in (0,\eta \Delta)$ with
$[r,s] \subseteq M_l(\lambda+\psi)$ for every $r,s \in R(\alpha)$
and $l \geq k$. But for $l < k$ if $r,s \in \ball{x}{\Delta}$ then by
the minimality of $k$ we have $\Delta \leq \psi w_n$ for every
$n \leq (1+\lambda)l$ so that by Lemma~\ref{shift},
$$
[r,s] \subseteq \ball{x}{\Delta} \subseteq M_l(\lambda+\psi).
$$
Hence for every $r,s \in R(\alpha) \cap \ball{x}{\Delta}$ we have
$
[r,s] \subseteq M_l(\lambda+\psi)
$
for any $l \geq 1$, so that $[r,s] \subseteq \T_{\lambda+\psi}$.
\end{proof}

We now let the Hilbert space $H$ equal $\R^d$ and write $S(H)$ for the unit sphere of $H$. Note that $(\setseq,\preceq) := ([0,1],\leq)$ is a dense, chain complete poset: for any $\lambda_1<\lambda_2$ there exists
$\lambda\in(\lambda_1,\lambda_2)$ and
every non-empty chain in $([0,1],\leq)$
has a supremum.

We quote \cite[Theorem~3.1]{DM} as Theorem~\ref{th.incr}. The assumptions
for Theorem~\ref{th.incr} are as
follows: $H$ is a real Hilbert space,
$(\setseq,\preceq)$ is a dense chain complete poset and
$(\T_{\lambda})_{\lambda\in\setseq}$ is a collection
of closed subsets of $H$ such that $\T_\lambda\subseteq \T_{\lambda'}$
whenever $\lambda\preceq\lambda'$.

We also use the following notation.
For a Lipschitz function $h\colon  H \rightarrow \R$ we write $D^h$ for the set of all pairs
$(x,e)\in H\times S(H)$ such that
the directional derivative $h'(x,e)$ exists and, for each $\lambda\in\setseq$,
we let $D^{h}_{\lambda}$  be the set of all $(x,e)\in D^h$ such that
$x\in \T_\lambda$.
If, in addition, $h\colon H\to\R$ is linear then we write $\|h\|$ for the operator norm of $h$.

\begin{theorem}
\label{th.incr}
Suppose $f_{0}\colon  H \rightarrow \R$ is a Lipschitz function, $\lambda_{0} \in\setseq$, $(x_{0},e_{0}) \in D^{f_{0}}_{\lambda_{0}}$, $\delta_{0},\mu,K > 0$ and $\lambda_{1} \in\setseq$ with $\lambda_{0} \prec \lambda_{1}$. Then there
exists a Lipschitz function $f\colon  H \rightarrow \R$
such that $f-f_{0}$ is linear with norm not greater than $\mu$ and a pair $(x,e) \in D^{f}_{\lambda}$, where
$\|x-x_{0}\| \leq \delta_{0}$ and $\lambda \in (\lambda_{0},\lambda_{1})$,
such that the directional derivative $f'(x,e) > 0$ is almost locally maximal in the following sense. For any
$\e > 0$ there exists $\delta_{\e} > 0$ and $\lambda_{\e} \in (\lambda,\lambda_{1})$ such that
whenever $(x',e') \in D^{f}_{\lambda_{\e}}$ satisfies
\begin{enumerate}
\item[(i)]
$\|x'-x\| \leq \delta_{\e}$, $f'(x',e') \geq f'(x,e)$ and
\item[(ii)]
for any $t\in\R$
\begin{equation}
\label{eqincr}
|(f(x'+te)-f(x'))-(f(x+te)-f(x))| \leq K \sqrt{f'(x',e')-f'(x,e)}|t|,
\end{equation}
\end{enumerate}
then we have $f'(x',e') < f'(x,e) + \e$.
\end{theorem}

We now quote \cite[Lemma~4.3]{DM}.
\begin{lemma}[Differentiability Lemma]\label{lemdiff}
Let $H$ be a real Hilbert space, $f:  H \rightarrow \R$ be a Lipschitz function and
$(x,e)\in H \times S(H)$ be such that the directional
derivative $f'(x,e)$ exists and is non-negative.
Suppose that there is a family of
sets $\{F_{\e}\subseteq H \mid\e>0\}$ such that
\begin{enumerate}
\item \label{cond-one}
whenever $\e,\eta > 0$ there exists
$\delta_*=\delta_*(\e,\eta) > 0$ such that for any $\delta \in(0, \delta_*)$ and $u_1,u_2,u_3$ in the closed unit ball of $H$, one can
find $u'_1,u'_2,u'_3$ with $\|u'_m-u_m\| \leq \eta$
and
\begin{equation*}
[x+\delta u'_1,x+\delta u'_3] \cup [x+\delta u'_3,x+\delta u'_2] \subseteq F_{\e},
\end{equation*}
\item\label{cond-two}
whenever $(x',e') \in F_\e \times S(H)$ is such that
the directional derivative $f'(x',e')$ exists,
$f'(x',e') \geq f'(x,e)$ and
\begin{align}\label{xx'}
&|(f(x'+te)-f(x'))-(f(x+te)-f(x))|\\
\notag &\leq
 25\sqrt{(f'(x',e')-f'(x,e))\mathrm{Lip}(f)}|t|
\end{align}
for every $t \in \R$ then
\begin{equation}\label{ineqless}
f'(x',e') < f'(x,e) + \e.
\end{equation}
\end{enumerate}
Then $f$ is \frechet{} differentiable at $x$ and its derivative
$f'(x)$ is given
by the formula
\begin{equation}\label{fgat}
f'(x)(h) = f'(x,e) \langle h,e \rangle
\end{equation}
for $h \in H$.
\end{lemma}

We now apply these results to our construction to obtain the following.

\begin{theorem}\label{thmmain}
There exists a compact subset $S \subseteq \R^d$ of Hausdorff dimension one with the universal differentiability property; moreover if $g \colon \R^d \to \R$ is Lipschitz, the set of points $x \in S$ such that $f$ is \frechet{} differentiable at $x$ is a dense subset of $S$.
\end{theorem}

\begin{proof}
We let
$$
S = \overline{\bigcup_{\q<1}\T_\q}.
$$
We first note 
$S\supseteq\T_0=\bigcap_{k=1}^\infty M_k(0)=R_1\ne\emptyset$
so that $S$ is not empty.
Also as $T_1$ is closed and $S \subseteq T_1 \subseteq O \subseteq \ball{0}{1}$, we conclude
$S$ is a compact set of Hausdorff dimension at most one. We shall prove it has the universal
differentiability property; it will follow, by Lemma~\ref{lem.low},
its Hausdorff dimension is equal to one.

Let $y \in S$, $\rho > 0$ and $g \colon \R^d \to \R$ be a Lipschitz function. We shall prove the existence of a point $x \in S$ of differentiability of $g$ with $\|x-y\| < \rho$.

We may assume $\mathrm{Lip}(g) > 0$. Let $H$ be the Hilbert space $\R^d$. We may pick $\q <\lambda_{1} := 1$ and $y' \in T_{\q}$ with $\|y'-y\| < \rho/3$.

Let $\lambda_0 \in (\q,1)$. By applying Lemma~\ref{lem.main} with $\eta < 1/2$ we can find distinct $r,s \in R \cap \ball{y'}{\rho/3}$ so that $[r,s]\subseteq\T_{\lambda_0}$.
Then, by Lebesgue's theorem, there exists $x_0 \in [r,s]$ such that $(x_{0},e_{0}) \in D^{g}_{\lambda_{0}}$, where
$$
e_0=\frac{r-s}{\|r-s\|}.
$$

Set $f_{0} = g$, $K = 25\sqrt{2\mathrm{Lip}(g)}$,
$\delta_{0} = \rho/3$ and $\mu = \mathrm{Lip}(g)$.

Let the Lipschitz function $f$,
the pair $(x,e)$, $\lambda\in(\lambda_0,\lambda_1)=(\lambda_0,1)$ and, for each $\e > 0$,
the numbers $\delta_{\e}>0$ and
$\lambda_{\e} \in (\lambda,1)$
be given by the conclusion of Theorem~\ref{th.incr}.
We verify
the conditions of  Lemma~\ref{lemdiff}
hold for the function $f\colon  \R^d \to \R$, the pair $(x,e) \in D^f_\lambda$ and the family of sets $\{F_{\e} \subseteq \R^d \mid \e > 0\}$ where
\begin{equation*}
F_{\e} =\T_{\lambda_{\e}} \cap \ball{x}{\delta_{\e}}.
\end{equation*}

We know from Theorem~\ref{th.incr} that the derivative $f'(x,e)$ exists and is non-negative.
To verify condition \eqref{cond-one} of Lemma~\ref{lemdiff},
for every $\e > 0$ and $\eta \in (0,1)$, we put
$\psi_\e=\lambda_\e-\lambda$ and define
$$
\delta_* =\frac{1}{2}
\min\Bigl\{
{\Delta_0\Bigl(\frac{\eta}{2},\psi_{\e}\Bigr)},\delta_{\e},1-\|x\|
\Bigr\},
$$
where 
$\Delta_0$ is given by
Lemma~\ref{lem.main}.
We see that $\delta\in(0,\delta_*)$ implies
$$
2\delta<\min\Bigl\{\Delta_0\Bigl(\frac{\eta}{2},\psi_{\e}\Bigr),1-\|x\|\Bigr\}.
$$
By Lemma~\ref{lem.main} we can find $\alpha\in(0,\eta/2 \cdot 2\delta)=(0,\eta\delta)$ such that for every $r,s \in R(\alpha) \cap B_{2\delta}(x)$ we have $[r,s] \subseteq T_{\lambda_{\e}}$. Using the definition of $R(\alpha)$ and $B_{2\delta}(x) \subseteq B_1(0)$ we can find
$x+\delta u_i'\in \ball{{x+\delta u_i}}{\alpha}$
such that
$$
[x+\delta u_1',x+\delta u_3']\cup[x+\delta u_3',x+\delta u_2'] \subseteq \T_{\lambda_\e}.
$$
Note then that since
$$
\|(x+\delta u_i')-(x+\delta u_i)\|<\alpha < \eta \delta,
$$
we have $\|u_i'-u_i\|<\eta$; also as
$\delta(1+\eta) < 2\delta_* \leq \delta_{\e}$
we have
$x+\delta u_i'\in\ball{x}{\delta_\e}$ for each $i=1,2,3$.
Thus
$$
[x+\delta u_1',x+\delta u_3']\cup[x+\delta u_3',x+\delta u_2'] \subseteq F_{\lambda_\e}.
$$
Condition \eqref{cond-two} of Lemma~\ref{lemdiff}
is immediate from the definition of
$F_{\e}$ and equation \eqref{eqincr} as $\mathrm{Lip}(f) \leq \mathrm{Lip}(g) + \mu = 2\mathrm{Lip}(g)$ so that
$25\sqrt{\mathrm{Lip}(f)} \leq K$.

Therefore, by Lemma~\ref{lemdiff} the function
$f$ is differentiable at $x$.
So too, therefore, is $g$ as $(g-f)$ is linear.
Finally, note that  $x \in \T_\lambda \subseteq S$ and
\begin{equation*}
\|x-y\| \leq \|x-x_{0}\| + \|x_{0}-y'\| + \|y'-y\| < \rho/3+\rho/3+\rho/3 = \rho.
\qedhere
\end{equation*}
\end{proof}
\bib
\end{document}